\documentclass[10pt]{amsart}
\usepackage{amsfonts,amssymb}
\def\qq{\mathbb{Q}}

\def\cc{\mathbb{C}}
\def\zz{\mathbb{Z}}

\def\pp{\mathbb{P}}
\def\hh{\mathcal{H}}

\def\tt{\mathcal{T}}
\def\ww{\mathcal{W}}

\def\difftials{H^{0}(X,\Omega_{X}^{1})}
\def\imtau{\mathrm{Im}\tau}

\theoremstyle{plain}
\newtheorem{theorem}{Theorem}[section]
\newtheorem{lemma}[theorem]{Lemma}
\newtheorem{proposition}[theorem]{Proposition}
\newtheorem{corollary}[theorem]{Corollary}

\theoremstyle{remark}
\newtheorem{remark}[theorem]{Remark}
\newtheorem{example}[theorem]{Example}

\theoremstyle{definition}
\newtheorem{definition}[theorem]{Definition}

\theoremstyle{definition}
\newtheorem{definition/proposition}[theorem]{Definition/Proposition}

\numberwithin{equation}{section}

\begin{document}

\title{Faltings' delta-invariant of a hyperelliptic Riemann surface}
\author{Robin de Jong}
\address{University of Amsterdam, The Netherlands}
\email{rdejong@science.uva.nl}

\begin{abstract}
In this note we give a closed formula for Faltings'
delta-invariant of a hyperelliptic Riemann surface.
\end{abstract}

\date{\today}

\maketitle

\section{Introduction}

In \cite{fa}, the delta-invariant is introduced for Riemann
surfaces, and it is suggested there to search for explicit
formulas for this invariant. In this note we give such an explicit
formula in the case of a hyperelliptic Riemann surface of arbitrary genus. We
note that \cite{fa} treats the case of elliptic curves, and that the case
of Riemann surfaces of genus 2 has been considered before in
\cite{bo}. 

In order to state our result, let us recall some notation and
earlier results. Let $X$ be a compact and connected Riemann
surface of genus $g>0$. Let $G$ be the Arakelov-Green function of
$X$ and let $\mu$ be the fundamental (1,1)-form on $X$
as defined in \cite{ar}, \cite{fa}.
Let $S(X)$ be the invariant defined by
\[ \log S(X) := - \int_X \log \|\vartheta\|(gP-Q) \cdot \mu(P) \, .
\] Here $\|\vartheta\|$ is the
function on $\mathrm{Pic}_{g-1}(X)$ defined as on \cite{fa},
p.~401, and $Q$ can be any point on $X$. The integral is well-defined and is
independent of the choice of the point $Q$.
In our paper \cite{jong}
we gave an explicit formula for the Arakelov-Green function of $X$.
\begin{theorem} \label{green}
Let $\ww$ be the classical divisor of Weierstrass points on $X$.
For $P,Q$ points on $X$, with $P$ not a Weierstrass point, we have
\[ G(P,Q)^g = S(X)^{1/g^2}
\cdot \frac{ \| \vartheta \|(gP-Q) }{ \prod_{W \in \ww} \|
\vartheta \|(gP-W)^{1/g^3}} \, . \] Here the product runs over the Weierstrass points of $X$,
counted with their weights. The formula is also valid if $P$ is a
Weierstrass point, provided that we take the leading
coefficients of a power series expansion about $P$ in both
numerator and denominator.
\end{theorem}
In the same paper, we also gave an explicit formula for the
delta-invariant of $X$. The delta-invariant is a fundamental
invariant of $X$, expressing the proportionality between two natural
metrics on the determinant of the Hodge bundle. For $P$ on $X$, not a Weierstrass point,
and $z$ a local coordinate about $P$, we put
\[ \|F_z\|(P) := \lim_{Q \to P} \frac{
\|\vartheta\|(gP-Q)}{|z(P)-z(Q)|^g} \, . \] Further we let
$W_z(\omega)(P)$ be the Wronskian at $P$ in $z$ of an orthonormal
basis $\{\omega_1,\ldots,\omega_g \}$ of the differentials
$\difftials$ provided with the standard hermitian inner product
$(\omega,\eta) \mapsto \frac{i}{2} \int_X \omega \wedge
\overline{\eta}$. We define an invariant $T(X)$ of $X$ by
\[ T(X)  := \| F_z \|(P)^{-(g+1)} \cdot \prod_{W
\in \ww} \| \vartheta \|(gP-W)^{(g-1)/g^3} \cdot |W_z(\omega)(P)|^2 \, ,
\] where again the product runs over the Weierstrass points of $X$,
counted with their weights. It can be checked that this does not
depend on the choice of $P$, nor on the choice of local coordinate
$z$ about $P$. A more intrinsic definition is possible, see
\cite{jong}, but the above formula will be convenient for us. We
remark that $T(X)$ does not involve an integral over $X$, contrary
to the invariant $S(X)$.
\begin{theorem} \label{robin}
For Faltings' delta-invariant $\delta(X)$ of $X$, the formula \[ \exp(\delta(X)/4) = S(X)^{-(g-1)/g^2} \cdot T(X) \]
holds.
\end{theorem}
In the present paper we make the invariant $T(X)$ explicit in the
case that $X$ is a hyperelliptic Riemann surface of genus $g \geq
2$. We relate it to a non-zero invariant
$\| \varphi_g\|(X)$ of $X$, the Petersson norm of the modular
discriminant associated to $X$, which we introduce in Section
\ref{modulardisc}. As we will see, for hyperelliptic 
Riemann surfaces this is a very
natural invariant to consider.
Unfortunately, it is not so clear how to extend its definition to
the general Riemann surface of genus $g$. \begin{definition} We
define by $G'$ the modified Arakelov-Green function
 \[ G'(P,Q) := S(X)^{-1/g^3} \cdot G(P,Q) \] on $X \times X$.
\end{definition}
We prove the following theorem dealing with $G'$ and $T(X)$.
Recall that the Weierstrass points of $X$ are just the ramification
points of a hyperelliptic map $X \to \pp^1$.
\begin{theorem} \label{main}
Let $W$ be a Weierstrass point of $X$. Let $n:={2g \choose g+1 }$.
Consider the product $ \prod_{W' \neq W} G'(W,W')$ running over
all Weierstrass points $W'$ different from $W$, ignoring their
weights. Then $ \prod_{W' \neq W} G'(W,W')$ is independent of the
choice of $W$ and the formula
\[  \prod_{W' \neq W} G'(W,W')^{(g-1)^2} = 2^{(g-1)^2}
\cdot \pi^{2g+2} \cdot T(X)^{\frac{g+1}{g}} \cdot \| \varphi_g
\|(X)^{\frac{1}{2n}}
\] holds.
\end{theorem}
The next theorem will be derived in a forthcoming article
\cite{jong4}. The result looks similar to the formula in Theorem
\ref{main}, but the proof is very different.
\begin{theorem} \label{second}
Let $m:={2g+2 \choose g}$. Then we have
\[ \prod_{(W,W')} G'(W,W')^{n (g-1)} =
\pi^{-2g(g+2)m} \cdot T(X)^{-(g+2)m} \cdot \| \varphi_g
\|(X)^{-\frac{3}{2}(g+1)} \, , \] the product running over all
ordered pairs of distinct Weierstrass points of $X$.
\end{theorem}
Combining the above two theorems yields a simple
closed formula for the invariant $T(X)$ in terms of
$\| \varphi_g \|(X)$.
\begin{theorem} Let $\|\Delta_g\|(X)$ be the modified
discriminant $\|\Delta_g\|(X) := 2^{-(4g+4)n} \cdot
\|\varphi_g\|(X)$. Then the formula \[ T(X) = (2 \pi)^{-2g} \cdot
\|\Delta_g\|(X)^{-\frac{3g-1}{8ng}} \] holds.
\end{theorem}
Combining this with Theorem \ref{robin} we obtain the following corollary.
\begin{corollary}
For Faltings' delta-invariant $\delta(X)$ of $X$, the formula
\[ \exp(\delta(X)/4) = (2 \pi)^{-2g} \cdot
S(X)^{-(g-1)/g^2} \cdot \|\Delta_g\|(X)^{-\frac{3g-1}{8ng}} \]
holds.
\end{corollary}
The significance of this result is that it makes the efficient
calculation of the delta-invariant possible for hyperelliptic
Riemann surfaces. We have given a demonstration of this in our
paper \cite{jong}.

We remark that in the case $g=2$, an explicit formula for the
delta-invariant has been given already by Bost \cite{bo}. Apart
from the Petersson norm of the modular discriminant, his formula
involves an invariant $\|H\|(X)$. This invariant has properties
similar to our $S(X)$.

The idea of the proof of Theorem \ref{main} is quite
straightforward: we start with the definition of the invariant
$T(X)$ and the formula for $G$ in Theorem \ref{green} and observe 
what happens if we let $P$ approach the Weierstrass point $W$ on $X$.
Thus, we have to perform a local study around $W$ of the function
$\prod_{W'} \|\vartheta\|(gP-W')$ and of the functions
$\|F_z\|(P)$ and $W_z(\omega)(P)$ for a suitable local coordinate
$z$. In Section \ref{localcoord} we find a suitable local
coordinate on an embedding of $X$ into its jacobian. In Section
\ref{expansions} we collect the local information that we need in
order to complete the proof in Section \ref{proofmainresult}. Some
preliminary work on this local information is carried out in the
Sections \ref{Schur} and \ref{sigma}. These two sections form the
technical heart of the paper.

\section{Modular discriminant} \label{modulardisc}

In this section we introduce the modular discriminant $\varphi_g$
and its Petersson norm $\|\varphi_g\|$. The modular discriminant
generalises the usual discriminant function $\Delta$ for elliptic
curves.

Let $g \geq 2$ be an integer and let $\hh_g$ be the Siegel upper
half-space of symmetric complex $g \times g$-matrices with
positive definite imaginary part. For $z \in \cc^g$ (viewed as a
column vector), a matrix $\tau \in \hh_g$ and $\eta,\eta' \in
\frac{1}{2} \zz^g$ we have the theta function with characteristic
$\eta=[{\eta' \atop \eta''}]$ given by
\[ \vartheta[\eta](z;\tau) := \sum_{n \in \zz^g} \exp( \pi i  {}^t (n+\eta')
\tau  (n+\eta') + 2\pi i  {}^t (n+\eta') (z+\eta'')) \, . \] For
any subset $S$ of $\{ 1,2,\ldots,2g+1\}$ we define a theta
characteristic $\eta_S$ as in \cite{mu}, Chapter IIIa: let
\[ \begin{array}{rcl} \eta_{2k-1} & = &
\left[ { {}^t( 0 \, , \, \ldots \, , \, 0 \, , \, \frac{1}{2} \, ,
\,  0 \, , \, \ldots \, , \, 0 ) \atop
{}^t (\frac{1}{2} \, , \, \ldots \, , \,  \frac{1}{2} \, , \,  0 \,
, \, 0 \, , \,
\ldots \, , \, 0 )}
\right] \, , \quad 1 \leq k \leq g+1 \, ,  \\ \eta_{2k} & = &
\left[ { {}^t( 0 \, , \, \ldots \, , \, 0 \, , \, \frac{1}{2} \, ,
\,  0 \, , \, 
\ldots \, , \, 0 ) \atop
{}^t (\frac{1}{2} \, , \, \ldots \, , \, \frac{1}{2} \, , \, \frac{1}{2} \, ,
\, 0 \, , \, \ldots \, , \, 0
)} \right] \, , \quad 1 \leq k \leq g \, ,
\end{array} \] where the non-zero entry in the top row occurs in
the $k$-th position. Then we put $\eta_S := \sum_{k \in S} \eta_k$
where the sum is taken modulo 1.
\begin{definition} \label{discrje} (\emph{Cf.} \cite{lock},
Section 3.) Let $\tt$ be the set of subsets of $\{1,2,
\ldots, 2g+1 \}$ of cardinality $g+1$. Write
$U:=\{1,3,\ldots,2g+1\}$ and let $\circ$ denote the symmetric
difference. The modular discriminant $\varphi_g$ is defined to be the function
\[ \varphi_g(\tau) := \prod_{T \in \tt}
\vartheta[\eta_{T \circ U}](0;\tau)^8   \] on $\hh_g$. The
function $\varphi_g$ is a modular form on $\Gamma_g(2) := \{ \gamma
\in \mathrm{Sp}(2g,\zz) | \gamma \equiv I_{2g} \bmod 2 \}$ of
weight $4r$ where $r:={2g+1 \choose g+1 }$.
\end{definition}
Consider an equation $y^2 = f(x)$ where $f \in \cc[X]$ is a monic
and separable polynomial of degree $2g+1$. Write $ f(x) =
\prod_{k=1}^{2g+1} (x-a_k) $ and denote by $D:=\prod_{k<l}
(a_k-a_l)^2$ the discriminant of $f$. Let $X$ be the hyperelliptic
Riemann surface of genus $g$ defined by $y^2=f(x)$. Then $X$
carries a basis of holomorphic differentials $ \mu_k := x^{k-1} dx
/2y $ where $k=1,\ldots,g$. Further, 
in \cite{mu}, Chapter IIIa, \S 5 it is
shown how, given an ordering of the roots of $f$, one can construct
a canonical symplectic basis of the homology of $X$. Throughout
this paper, we will always work with such a canonical basis of
homology, \emph{i.e.}, a certain ordering of the roots of a
hyperelliptic equation will always be taken for granted.

Let $(\mu|\mu')$ be the period matrix
of the differentials $\mu_k$ with respect to a chosen canonical
basis of homology. Then $\mu$ is invertible, and we put $\tau
:= \mu^{-1} \mu'$.
\begin{proposition} \label{disc}
We have the formula
\[ D^n = \pi^{4gr} (\det \mu)^{-4r} \varphi_g(\tau)   \]
relating the discriminant $D$ of the polynomial $f$ 
to the value $\varphi_g(\tau)$ of the modular discriminant.
\end{proposition}
\begin{proof}
See \cite{lock}, Proposition 3.2. 
\end{proof}
\begin{definition}
Let $X$ be a hyperelliptic Riemann surface of genus $g \geq 2$
and let $\tau$ be a period matrix for $X$ formed on a canonical
symplectic basis, given by an ordering of the roots of an 
equation $y^2=f(x)$ for $X$. Then we
write $\|\varphi_g\|(\tau)$ for the Petersson norm $(\det
\imtau)^{2r} \cdot | \varphi_g(\tau)|$ of $\varphi_g(\tau)$. 
This does not depend on the choice of
$\tau$ and hence it defines an invariant $\|\varphi_g\|(X)$ of $X$.
\end{definition}
It follows from Proposition \ref{disc} that the invariant
$\|\varphi_g\|(X)$ is non-zero.

\section{Local coordinate} \label{localcoord}

For our local computations on our hyperelliptic Riemann surface 
we need a convenient local coordinate.
We find one by embedding the Riemann surface into its
jacobian and by taking one of the euclidean coordinates.

Let $X$ be a hyperelliptic Riemann surface of genus $g \geq 2$,
let $y^2=f(x)$ with $f$ monic of degree $2g+1$ be an equation for
$X$, let $\mu_k$ be the differential given by $\mu_k =
x^{k-1}dx/2y$ for $k=1,\ldots,g$, and let $(\mu|\mu')$ be their
period matrix formed on a canonical basis of homology. Let $L$ be
the lattice in $\cc^g$ generated by the columns of $(\mu|\mu')$.
We have an embedding $\iota : X \hookrightarrow \cc^g/L$ given by
integration $P \mapsto \int_{\infty}^P (\mu_1,\ldots,\mu_g)$. We
want to express the coordinates $z_1,\ldots,z_g$, restricted to
$\iota(X)$, in terms of a local coordinate about
$0=\iota(\infty)$. This is established by the following lemma. In
general, we denote by $O(w_1,\ldots,w_s;d)$ a Laurent series in
the variables $w_1,\ldots,w_s$ all of whose terms have total
degree at least $d$. We owe the argument to \cite{on}.
\begin{lemma} \label{coordinates}
The coordinate $z_g$ is a local coordinate about 0 on $\iota(X)$,
and we have
\[  z_k = \frac{1}{2(g-k)+1}z_g^{2(g-k)+1} + O(z_g;2(g-k)+2)  \] on $\iota(X)$ for $k=1,\ldots,g$.
\end{lemma}
\begin{proof} We can choose a local coordinate $t$ about $\infty$
on $X$ such that $x=t^{-2}$ and $y=-t^{-(2g+1)}+O(t;-2g)$. For $P
\in X$ in a small enough neighbourhood of $\infty$ on $X$ and for
a suitable integration path on $X$ we then have
\[ \begin{aligned} z_k(P) &= \int_{\infty}^P \frac{x^{k-1} dx}{2y} =
\int_0^{t(P)} \frac{ t^{-2(k-1)} \cdot (-2t^{-3} dt) }{-2t^{-(2g+1)}
+ O(t;-2g) } \\ &= \int_0^{t(P)} \left( t^{2(g-k)} +
O(t; 2(g-k)+1) \right) dt \\ &= \frac{1}{2(g-k)+1}
t(P)^{2(g-k)+1} + O(t(P); 2(g-k)+2 ) \, .
\end{aligned} \] By taking $k=g$ we find $z_g = t + O(t;
2)$ and for $k=1,\ldots,g-1$ then
\[ z_k = \frac{1}{2(g-k)+1} z_g^{2(g-k)+1} + O(z_g;
2(g-k)+2) \, , \] which is what we wanted.
\end{proof}

\section{Schur polynomials} \label{Schur}

In this section we assemble some facts on Schur polynomials. We
will need these facts at various places in the next sections. Fix
a positive integer $g$. Consider the ring of symmetric polynomials
with integer coefficients in the variables $x_1,\ldots,x_g$.
Let $e_r$ be the elementary symmetric functions defined by means of the
generating function $E(t) = \sum_{r \geq 0} e_r t^r = \prod_{k=1}^g 
(1+x_k t)$. 
\begin{definition} \label{Sg} Let $d$ be a positive integer
and let $\pi=\{\pi_1,\ldots,\pi_h \}$ with $\pi_1 \geq \ldots \geq
\pi_h$ be a partition of $d$. The
Schur polynomial associated to $\pi$ is the polynomial
\[ S_\pi := \det( e_{\pi'_k-k+l} )_{1 \leq k,l \leq h} \, , \]
where $h$ is the length of the partition $\pi$, and 
where $\pi'$ is the conjugate partition of $\pi$ given by 
$\pi'_k = \# \{ l \, : \, \pi_l \geq k \}$, \emph{i.e.},
the partition obtained by
switching the associated Young diagram around its diagonal. 
The polynomial
$S_\pi$ is symmetric and has total degree $d$. We denote by $S_g$
the Schur polynomial in $g$ variables associated to the partition
$\pi = \{g,g-1,\ldots,2,1 \}$. Thus, the formula
\[  S_g = \det \left( e_{g-2k+l+1} \right)_{1 \leq k,l \leq g}   \]
holds, and the polynomial $S_g$ has total degree $g(g+1)/2$.
\end{definition}
Let $p_r$ be the elementary Newton functions (power
sums) given by the generating function $P(t) = \sum_{r \geq 1}
p_r t^{r-1} = \sum_{k \geq 1} x_k/(1-x_k t)$. The following proposition 
is then a special case of Theorem 4.1 of
\cite{buch}.
\begin{proposition} \label{unique}
The Schur polynomial $S_g$ can be expressed as a polynomial in the
$g$ functions $p_1,p_3,\ldots,p_{2g-1} $ only. This polynomial is
unique.
\end{proposition}
\begin{definition} \label{esg}
We define $s_g$ to be the unique polynomial in
$g$ variables given by Proposition \ref{unique}.
\end{definition}
The next proposition is a special case of Theorem 6.2 of
\cite{buch}.
\begin{proposition} \label{Schurproperty}
Let $s(x_1,\ldots,x_g) \in \cc[x_1,\ldots,x_g]$ be a polynomial in
$g$ variables such that for any set of $g$ complex numbers
$w_1,\ldots,w_g$, the polynomial $
s(z_1-w,z_2-w^3,\ldots,z_g-w^{2g-1}) $ in $w$ either has exactly
$g$ roots $w_1,\ldots,w_g$, or vanishes identically, if we give $z$ 
the value
$z=(p_1(w_1,\ldots,w_g),p_3(w_1,\ldots,w_g), \ldots,
p_{2g-1}(w_1,\ldots,w_g))$. Then $s$ is equal to the polynomial $s_g$
up to a constant factor.
\end{proposition}
\begin{definition} \label{sigmag}
We define $\sigma_g$ to be the polynomial in
$g$ variables given by the equation
\[ \sigma_g(z_1,\ldots,z_g) = s_g(z_g,3z_{g-1},\ldots,(2g-1)z_1) \,
. \]
\end{definition}
The following proposition is then the result of a simple calculation.
\begin{proposition} \label{Hankel}
Up to a sign, the homogeneous part of least total degree of
$\sigma_g$ is equal to the Hankel determinant
\[ H(z) = \det \left( \begin{array}{cccc} z_1 & z_2 & \cdots &
z_{(g+1)/2} \\ z_2 & z_3 & \cdots & z_{(g+3)/2} \\ \vdots & \vdots
& \ddots & \vdots \\ z_{(g+1)/2} & z_{(g+3)/2} & \cdots & z_g
\end{array} \right) \] if $g$ is odd, or
\[ H(z) = \det \left( \begin{array}{cccc} z_1 & z_2 & \cdots &
z_{g/2} \\ z_2 & z_3 & \cdots & z_{(g+2)/2} \\ \vdots & \vdots &
\ddots & \vdots \\ z_{g/2} & z_{(g+2)/2} & \cdots & z_{g-1}
\end{array} \right) \] if $g$ is even.
\end{proposition}
We conclude with some more general facts. These can all be found
for example in Appendix A to \cite{fulton}.
\begin{proposition} \label{at111}
Let $\pi = \{ \pi_1,\ldots,\pi_h \}$ with $\pi_1 \geq \ldots \geq
\pi_h $ be a partition. Then the
formula \[ S_\pi(1,\ldots,1) = \prod_{k<l} \frac{\pi_k - \pi_l + l
-k}{l - k} \] holds. 
In particular, $S_g(1,\ldots,1)=2^{g(g-1)/2}$.
\end{proposition}
\begin{definition} Let $\mathbf{i}=(i_1,\ldots,i_d)$ be a $d$-tuple
of non-negative integers. The $\mathbf{i}$-th generalised Newton
function $p^{(\mathbf{i})}$ is defined to be the polynomial
\[ p^{(\mathbf{i})} := p_1^{i_1} \cdot p_2^{i_2} \cdot \ldots \cdot
p_d^{i_d} \, , \] where the $p_r$ are the elementary Newton
functions.
\end{definition}
\begin{proposition}
The set of generalised Newton functions $p^{(\mathbf{i})}$, where
$\mathbf{i}$ runs through the $d$-tuples
$\mathbf{i}=(i_1,\ldots,i_d)$ of non-negative integers with $\sum
\alpha i_\alpha =d$, forms a basis of the $\qq$-vector space of
symmetric polynomials with rational coefficients 
of total degree $d$.
\end{proposition}
\begin{proposition} \label{expandSchur}
For a partition $\pi$ of $d$ and a $d$-tuple
$\mathbf{i}=(i_1,\ldots,i_d)$, denote by $\omega_\pi(\mathbf{i})$
the coefficient of the monomial $x_1^{\pi_1} \cdot \ldots \cdot
x_d^{\pi_d}$ in $p^{(\mathbf{i})}$. Then the polynomial $S_\pi$
can be expanded on the basis $\{ p^{(\mathbf{i})} \}$ of
generalised Newton functions of total degree $d$ as $ S_\pi =
\sum_\mathbf{i} \frac{1}{z(\mathbf{i})} \cdot
\omega_\pi(\mathbf{i}) \cdot p^{(\mathbf{i})} $. Here
$z(\mathbf{i}) = i_1!1^{i_1} \cdot i_2! 2^{i_2} \cdot \ldots \cdot
i_d! d^{i_d} $.
\end{proposition}

\section{Sigma function} \label{sigma}

We consider again hyperelliptic Riemann surfaces of genus $g \geq
2$, defined by equations $y^2 = f(x)$ with $f$ monic and separable
of degree $2g+1$. We write $ f(x) = x^{2g+1} + \lambda_1 x^{2g} +
\cdots + \lambda_{2g}x + \lambda_{2g+1} $  and denote by $\lambda$
the vector of coefficients $(\lambda_1,\ldots,\lambda_{2g+1})$. In
this section we study the sigma function $\sigma(z;\lambda)$ with
argument $z \in \cc^g$ and parameter $\lambda$. This is a modified
theta function, studied extensively in the nineteenth century. 
Klein observed that the sigma function serves very well to
study the function theory of hyperelliptic Riemann surfaces. For
us it will be a convenient technical tool for obtaining the local
expansions that we need. We will give the definition of the sigma
function, as well as its power series expansion in $z,\lambda$. For more
details we refer to the \emph{Enzyklop\"adie der mathematischen 
Wissenschaften}, Band II, Teil 2, Kapitel 7.XII.
A modern reference is \cite{bel3}, where one also finds
applications of the sigma function in the theory of the
Korteweg-de Vries differential equation.

As before, let $\mu_k$ be the holomorphic differential given by
$\mu_k = x^{k-1}dx/2y$ for $k=1,\ldots,g$, and let $(\mu|\mu')$ be
their period matrix formed on a canonical basis of homology. Let
$L$ be the lattice in $\cc^g$ generated by the columns of
$(\mu|\mu')$. By the theorem of Abel-Jacobi we have a bijective
map $ \mathrm{Pic}_{g-1}(X) {\buildrel \sim \over \longrightarrow}
\cc^g/L$ given by $ \sum_k m_k P_k \longmapsto \sum_k m_k
\int_{\infty}^{P_k} (\mu_1,\ldots,\mu_g) $. Denote by $\Theta$ the
image of the theta divisor of classes of effective divisors of
degree $g-1$, and let $q:\cc^g \to \cc^g/L$ be the projection map.
Let $\tau = \mu^{-1}\mu'$. By a fundamental theorem of Riemann,
there exists a unique theta-characteristic $\delta$ such that
$\vartheta[\delta](z;\tau)$ vanishes to order one precisely along
$q^{-1}(\Theta)$. 
\begin{definition} Let $\nu$ be the matrix of $A$-periods of the
differentials of the second kind $ \nu_k := \frac{1}{4y}
\sum_{l=k}^{2g-k} (l+1-k)\lambda_{l+k+1} x^k dx $ for
$k=1,\ldots,g$. These differentials have a second order pole at
$\infty$ and no other poles. The sigma function is then the
function
\[  \sigma(z;\lambda) :=
\exp(-\frac{1}{2}z\nu \mu^{-1} {}^t z) \cdot
\vartheta[\delta](\mu^{-1}z;\tau) \, . \]
\end{definition}
Using some of the facts on Schur polynomials from the previous
section, we can give the power series expansion of
$\sigma(z;\lambda)$. The result is probably well-known to specialists, although
we couldn't find an explicit reference in the literature. For the formulation
and the proof we were inspired by \cite{on}, as well as by a private
communication with the author.
For the special case $g=2$, a somewhat
stronger version of the result has
been obtained by Grant, see \cite{gr}, Theorem 2.11.
\begin{proposition} \label{SigmaTaylor}
The power series expansion of $\sigma(z;\lambda)$ about $z=0$ is
of the form
\[ \sigma(z;\lambda) = \gamma \cdot
\sigma_g(z) + O(\lambda) \, ,
\] where $\sigma_g$ is the polynomial given by Definition \ref{sigmag} and
where the symbol $O(\lambda)$ denotes a power series in
$z,\lambda$ in which each term contains a $\lambda_k$ raised to a
positive integral power. The constant $\gamma$ satisfies the
formula
\[ \gamma^{8n} =  \pi^{4g(r-n)} (\det \mu)^{-4(r-n)} \varphi_g(\tau) \, . \]
If we assign the variable $z_k$ a weight $2(g-k)+1$, and the
variable $\lambda_k$ a weight $-2k$, then the power series
expansion in $z,\lambda$ of $\sigma(z;\lambda)$ is homogeneous of
weight $g(g+1)/2$.
\end{proposition}
\begin{proof}  First of all, the
homogeneity of the power series expansion in $z,\lambda$ with
respect to the assigned weights follows from an explicit formula
for $\sigma(z;\lambda)$ given in \cite{bel2}. This homogeneity
is also mentioned there, \emph{cf.} the concluding remarks after
Corollary 1. Write $\sigma(z;\lambda) = \sigma_0(z) + O(\lambda)$
where $O(\lambda)$ denotes a power series in $z,\lambda$ in which
each term contains a $\lambda_k$ raised to a positive integral
power. Because of the homogeneity, the series $\sigma_0(z)$ is
necessarily a polynomial in the variables $z_1,\ldots,z_g$. By the
Riemann vanishing theorem, there is a dense open subset $U \subset
\cc^{2g+1}$ such that for any $\lambda \in U$, the function
$\sigma(z;\lambda)$ satisfies the following property: for any set
of $g$ points $P_1,\ldots,P_g$ on the hyperelliptic Riemann
surface $X=X_\lambda$ corresponding to $\lambda$, the function
$\sigma(z - \int_\infty^P(\mu_1,\ldots,\mu_g);\lambda)$ in $P$ on
$X$ either has exactly $g$ roots $P_1,\ldots,P_g$, or vanishes
identically, when we give the argument $z$ the value $z = \sum_k
\int_\infty^{P_k}(\mu_1,\ldots,\mu_g)$.
In the limit $\lambda \to
0$ we find then, as in the proof of Lemma \ref{coordinates}, that
for any set of $g$ complex numbers $w_1,\ldots,w_g$ the polynomial
\[ \sigma_0 \left( \frac{1}{2g-1} \left(z_g - w^{2g-1} \right),
\frac{1}{2g-3} \left( z_{g-1} - w^{2g-3} \right), \ldots,
\frac{1}{3} \left(z_2 - w^3 \right), z_1 - w \right) \] in $w$
either has exactly $g$ roots $w_1,\ldots,w_g$, or vanishes
identically, for the value
$z=(p_1(w_1,\ldots,w_g),p_3(w_1,\ldots,w_g), \ldots,
p_{2g-1}(w_1,\ldots,w_g))$. By Proposition \ref{Schurproperty},
the polynomial $\sigma_0$ is equal to the polynomial $\sigma_g$ up to a
constant factor $\gamma$. As to this constant $\gamma$, we find in
\cite{ba}, Section IX a calculation of a constant $\gamma'$ such
that $ \sigma(z;\lambda) = \gamma' \cdot H(z) + O(z;\lfloor
(g+3)/2 \rfloor) $, where $H(z)$ is the Hankel determinant from
Proposition \ref{Hankel} and where now we consider the power
series expansion only with respect to the variables
$z_1,\ldots,z_g$ and with respect to their usual weight
$\deg(z_k)=1$. By Proposition \ref{Hankel}, this $\gamma'$ is equal to
our constant $\gamma$, up to a sign. We just quote the result of
Baker's computation: \[ \gamma^4 =
 \vartheta(0;\tau)^4 \cdot \prod_{ {k<l \atop k,l \in U} }
(a_k - a_l)^2 /( \ell_1 \ell_3 \cdots \ell_{2g+1} ) \, , \,
\textrm{where} \, \ell_r
:= -i \cdot \prod_{ {k \in U \atop k \neq r} } (a_k - a_r) /
\prod_{k \notin U} (a_k - a_r) \, . \] Thomae's formula (\emph{cf.}
\cite{mu}, Chapter IIIa, \S 8) says that \[
\vartheta(0;\tau)^8 = (\det \mu)^4 \pi^{-4g} \prod_{ {k < l \atop
k,l \in U }} (a_k - a_l)^2 \prod_{ {k<l \atop k,l \notin U}} (a_k
- a_l)^2 \, . \] 
Combining, we obtain $ \gamma^8 = D \cdot \pi^{-4g}
\cdot (\det \mu)^4 $. The formula for $\gamma$ that we gave then
follows from Proposition \ref{disc}.
\end{proof}
\begin{example}
By way of illustration, we have computed $\sigma_g$ for small $g$:
\\
\smallskip
\medskip\noindent
\vbox{
\bigskip\centerline{\def\quad{\hskip 0.6em\relax}
\def\quod{\hskip 0.5em\relax }
\vbox{\offinterlineskip \hrule
\halign{&\vrule#&\strut\quod\hfil#\quad\cr
height2pt&\omit&&\omit&\cr & $g$ &&$\sigma_g$&\cr
height2pt&\omit&&\omit&\cr \noalign{\hrule}
height2pt&\omit&&\omit&\cr &1 && $z_1$ &\cr
height2pt&\omit&&\omit&\cr &2 && $-z_1 + \frac{1}{3}z_2^3$ &\cr
height2pt&\omit&&\omit&\cr &3 && $z_1z_3 - z_2^2 - \frac{1}{3} z_2
z_3^3 + \frac{1}{45}
            z_3^6 $ &\cr
height2pt&\omit&&\omit&\cr &4 && $ z_1z_3 - z_2^2 -z_3^2z_4 +
z_2z_3z_4^2 - \frac{1}{3} z_1z_4^3 +
    \frac{1}{15}z_2 z_4^5 - \frac{1}{105} z_3 z_4^7 + \frac{1}{4725}z_4^{10}$
    &\cr
height2pt&\omit&&\omit&\cr }  \hrule } }}
\end{example}
\begin{remark} As can be seen from Proposition \ref{Hankel},
the homogeneous part of least total degree (with respect to the
usual weight $\deg(z_k)=1$) of $\sigma_g(z)$ has degree $\lfloor
(g+1)/2 \rfloor$. Hence, by a fundamental theorem of Riemann, the
theta-characteristic $\delta$ gives rise to a linear system of
dimension $\lfloor (g-1)/2 \rfloor$ on $X$. 
\end{remark}

\section{Leading coefficients}
\label{expansions}

In this section we calculate the leading coefficients of the power
series expansions in $z_g$ of the holomorphic functions
$\vartheta[\delta](g \mu^{-1}z;\tau)|_{\iota(X)}$ and
$W_{z_g}(\mu)$, the Wronskian in $z_g$ of the basis
$\{\mu_1,\ldots,\mu_g \}$.
\begin{proposition} \label{SigmaHolo}
The leading coefficient of the power series expansion of $\sigma(g
z;\lambda) |_{\iota(X)} $, and hence of $\vartheta[\delta](g
\mu^{-1}z;\tau)|_{\iota(X)}$, is equal to $ \gamma \cdot 2^{g(g-1)/2} $.
\end{proposition}
\begin{proof}
By Lemma \ref{coordinates} and Proposition \ref{SigmaTaylor}, the
power series expansion of $\sigma(g z;\lambda) |_{\iota(X)} $ has
the form
\[ \begin{aligned} \sigma(g z;\lambda)  
|_{\iota(X)} &=  \gamma \cdot \sigma_g \left(
\frac{g}{2g-1}z_g^{2g-1},\frac{g}{2g-3}z_g^{2g-3},\ldots,\frac{g}{3}z_g^3,gz_g
\right) \\ &+ O(z_g; g(g+1)/2+1) \, . \end{aligned} \] Hence we need
to calculate $ \sigma_g
\left(\frac{g}{2g-1},\frac{g}{2g-3},\ldots,\frac{g}{3},g \right)$.
By Definition \ref{sigmag} this is equal to $ s_g(g,g,\ldots,g) $. But by
Proposition \ref{unique} and Definition \ref{esg} we have
$s_g(g,g,\ldots,g)=S_g(1,1,\ldots,1)$, and by Proposition
\ref{at111} we have $S_g(1,\ldots,1)=2^{g(g-1)/2}$. The
proposition follows.
\end{proof}
\begin{proposition} \label{WronsHolo} The leading coefficient of
the power series expansion of the Wronskian $ W_{z_g}(\mu) $
is equal to $ \pm 2^{g(g-1)/2}$.
\end{proposition}
\begin{proof} Expanding the Wronskian yields
\[ \begin{aligned}
    W_{z_g}(\mu)  &= \det \left( \frac{1}{(k-1)!}
     \frac{d^k z_l}{dz_g^l} \right)_{1\leq k,l \leq g} = \\
     &= \left( \begin{array}{ccccc}
        z_g^{2g-2} & z_g^{2g-4} & \cdots & z_g^2 & 1 \\
        (2g-2)z_g^{2g-3} & (2g-4)z_g^{2g-5} & \cdots & 2z_g & 0 \\
        \vdots & \vdots & \ddots & \vdots & \vdots \\
        {2g-2 \choose g-1}z_g^g & {2g-4 \choose g-1}z_g^{g-2} &
        \cdots & 0 & 0 \end{array} \right) + O(z_g; g(g-1)/2+1)
    \, . \end{aligned} \]
Let $A$ be the matrix of binomial coefficients $ A := \left( {2g-2k
\choose g-l} \right)_{1 \leq k,l \leq g-1} $. From the expansion
it follows that the required leading coefficient is equal to $\det A$.
We will compute this number. First of all note that
\[ \det A = \frac{ (2g-2)!(2g-4)! \cdots 2! }{(g-1)!(g-2)!
\cdots 1!} \det \left( \frac{1}{(g-2k+l)!} \right)_{1\leq k,l \leq
g-1} \, ,
\] where we define $1/n! := 0$ for $n<0$. Now let $d=g(g-1)/2$ and
consider the ring of symmetric polynomials with integer
coefficients in $g-1$ variables. It is well known that for the
elementary symmetric functions $e_r$ we have an expansion
\[ e_r = \frac{1}{r!} \det \left( \begin{array}{ccccc}
        p_1  &  1  &  0  &  \cdots  &  0  \\
        p_2  &  p_1 & 2  &  \cdots  & 0  \\
        \cdots & \cdots  &  \cdots &  \cdots  & \cdots    \\
        p_{r-1} & p_{r-2} & p_{r-3} & \cdots & r-1 \\
        p_r  & p_{r-1} & p_{r-2} & \cdots & p_1 \end{array} \right) \, ,
\] with $p_r$ the elementary Newton functions. From Definition
\ref{Sg} and this expansion it follows that $\det(1/(g-2k+l)!)$ is
the coefficient of $p_1^d$ in the expansion of $S_{g-1}$ with
respect to the basis of generalised Newton functions. By
Proposition \ref{expandSchur}, this coefficient is equal to $
\omega_{g-1}(d)/d!$, where $\omega_{g-1}(d)$ is
the coefficient of $x_1^{g-1}   x_2^{g-2} \cdots x_{g-1}^2 x_g$ in
$ p_1^d$. It immediately follows
that $\det(1/(g-2k+l)!) = 1/(g-1)!(g-2)! \cdots 1!  $.
Combining one finds $\det A = 2^{g(g-1)/2}$.
\end{proof}

\section{Proof of Theorem \ref{main}} \label{proofmainresult}

Now we are ready to prove Theorem \ref{main}. Let $X$ be a
hyperelliptic Riemann surface of genus $g \geq 2$, and let $W$ be
one of its Weierstrass points.
\begin{proof}[Proof of Theorem \ref{main}]
Fix a hyperelliptic equation $y^2=f(x)$ 
for $X$ with $f$ monic and separable of degree $2g+1$ that puts $W$ at infinity.
Choose a canonical basis of the homology of $X$, and form the
period matrix $(\mu|\mu')$ of the differentials $x^{k-1} dx/2y$
for $k=1,\ldots,g$ on this basis. Let $L$ be the lattice in
$\cc^g$ generated by the columns of $(\mu|\mu')$, and embed $X$
into $\cc^g/L$ with base point $W$ as in Section \ref{localcoord}.
We have the standard euclidean coordinates $z_1,\ldots,z_g$ on
$\cc^g/L$ and according to Lemma \ref{coordinates} we have that
$z_g$ is a local coordinate about $W$ on $X$. The weight $w$ of
$W$ is given by $w=g(g-1)/2$. Consider then the following
quantities:
\[ \begin{aligned}
A(W') &:= \lim_{Q \to W} \frac{ \| \vartheta \|(gQ-W') }{ |z_g|^g}
\quad \textrm{for Weierstrass points} \quad  W' \neq W \, ; \\
A(W) &:= \lim_{Q \to W} \frac{ \| \vartheta \|(gQ-W) }{
|z_g|^{w+g} } = \lim_{Q \to W} \frac{ \| F_{z_g} \|(Q) }{ |z_g|^w
}  \, ;
\\ B(W) &:= \lim_{Q \to W} \frac{ | W_{z_g}(\omega)(Q) |}{
|z_g|^w} \, , \end{aligned} \] where $W_{z_g}(\omega)$ is the
Wronskian in $z_g$ of an orthonormal basis $\{
\omega_1,\ldots,\omega_g \}$ of $\difftials$. We have by Theorem
\ref{green}
\[ G'(W,W')^g = \frac{A(W')}{
 \prod_{W''} A(W'')^{w/g^3}} \quad \textrm{for
Weierstrass points} \quad  W' \neq W \, , \] hence
\[ \prod_{W' \neq W} G'(W,W')^g = \frac{1}{A(W)} \cdot \left( \prod_{W'}
A(W') \right)^{\frac{g+1}{2g^2}} \, . \] Further we have by the
definition of $T(X)$, letting $P$ approach $W$, \[ T(X) =
A(W)^{-(g+1)} \cdot \left( \prod_{W'} A(W')
\right)^{\frac{w(g-1)}{g^3}} \cdot B(W)^2 \, . \] Eliminating the
factor $\prod_{W'} A(W')$ yields
\[ \prod_{W'\neq W} G'(W,W')^{(g-1)^2} = A(W)^4 \cdot
B(W)^{-\frac{2g+2}{g}} \cdot T(X)^{\frac{g+1}{g}} \, . \] Now we
use the results obtained in Section \ref{expansions}. Let $\tau =
\mu^{-1}\mu'$. A simple calculation gives that $ A(W)$ is $(\det
\imtau)^{1/4}$ times the absolute value of the leading coefficient
of the power series expansion of $ \vartheta[\delta](g
\mu^{-1}z;\tau)|_{\iota(X)}$ in $z_g$. Hence by Propositions
\ref{SigmaTaylor} and \ref{SigmaHolo} we have
\[ A(W) =  2^{g(g-1)/2} \cdot
\pi^{g \frac{r-n}{2n}} \cdot (\det \imtau)^{1/4} \cdot |\det
\mu|^{-\frac{r-n}{2n}} \cdot
 | \varphi_g(\tau) |^{\frac{1}{8n}} \, . \]
Next let $\| \cdot \|$ be the metric on $\wedge^g \difftials$
derived from the hermitian inner product $(\omega,\eta) \mapsto
\frac{i}{2} \int_X \omega \wedge \overline{\eta}$ on $\difftials$.
We then have $\| \mu_1 \wedge \ldots
\wedge \mu_g \|^2 =   (\det \imtau) \cdot |\det \mu|^2$ by Riemann's
second bilinear relations. This
gives that $|W_{z_g}(\omega)| = |W_{z_g}(\mu)| \cdot (\det
\imtau)^{-1/2} \cdot |\det \mu|^{-1}$. From Proposition
\ref{WronsHolo} we derive then
\[ B(W) = 2^{g(g-1)/2} \cdot (\det \imtau)^{-1/2} \cdot |\det \mu|^{-1} \, . \]
Plugging in our results for $A(W)$ and $B(W)$ finally gives the
theorem.
\end{proof}
\begin{remark}
The fact that the product from Theorem \ref{main} is independent
of the choice of the Weierstrass point $W$ follows \emph{a
fortiori} from the computations in the above proof. It would be
interesting to have an \emph{a priori} reason for this
independence.
\end{remark}
\begin{remark}
We have not been able to find in general a formula for $G'(W,W')$
with $W,W'$ just two Weierstrass points. In the case $g=2$ it can
be shown that \[ G'(W,W')^2 = 2^{1/4} \cdot
\|\varphi_2\|(X)^{-3/64} \cdot \prod_{W'' \neq W,W'}
\|\vartheta\|(W-W'+W'') \, . \] This formula should be compared
with the explicit formula for $G(W,W')$ given in \cite{bo},
Proposition 4. We guess that in general we have \[ (??) \quad
G'(W,W')^g = A(X) \cdot \prod_{ \{W_1,\ldots,W_{g-1} \} \, , \atop
W,W' \notin \{W_1,\ldots,W_{g-1} \} }
\|\vartheta\|(W-W'+W_1+\cdots+W_{g-1}) \, , \] with $A(X)$ some
invariant of $X$. Such a result is consistent with Theorems
\ref{main} and \ref{second}.
\end{remark}

\subsection*{Acknowledgments} The author wishes to thank his thesis
advisor Gerard van der Geer for his encouragement and helpful
remarks. He also kindly thanks Professor Yoshihiro \^Onishi for
his remarks pertaining to Proposition \ref{SigmaTaylor}.

\end{document}